\documentclass[11pt,reqno]{amsart}
\usepackage{amsmath,amssymb,amsthm}
\usepackage[a4paper, margin=2cm]{geometry}
\usepackage[colorlinks=true,linkcolor=blue,citecolor=blue]{hyperref}
\numberwithin{equation}{section}

\theoremstyle{plain}
\newtheorem{theorem}{Theorem}[section]
\newtheorem{proposition}[theorem]{Proposition}
\newtheorem{lemma}[theorem]{Lemma}
\newtheorem{corollary}[theorem]{Corollary}

\newcommand{\Q}{\mathbb{Q}}
\newcommand{\Z}{\mathbb{Z}}
\newcommand{\R}{\mathbb{R}}
\newcommand{\OK}{\mathcal{O}_K}
\newcommand{\Ov}{\mathcal{O}_K^{\vee}}
\newcommand{\dif}{\mathfrak{d}}
\newcommand{\qq}{\mathfrak{q}}
\newcommand{\PP}{\mathfrak{P}}
\newcommand{\Tr}{\operatorname{Tr}}
\newcommand{\Nm}{\operatorname{N}}
\newcommand{\disc}{\operatorname{disc}}
\newcommand{\covol}{\operatorname{covol}}
\newcommand{\vol}{\operatorname{vol}}
\newcommand{\rank}{\operatorname{rank}}
\newcommand{\tp}{\succ 0}
\newcommand{\pr}[1]{\textup{(\ref{#1})}}

\begin{document}

\title[Universal quadratic forms have unbounded rank in every degree]
{Universal quadratic forms have unbounded rank\\ in every degree}

\author{Pavlo Yatsyna}
\address{Charles University, Faculty of Mathematics and Physics, Department of Algebra, Sokolovská 49/83, 186 75 Praha 8, Czech Republic}
\email{p.yatsyna@matfyz.cuni.cz}


\begin{abstract}
We prove that for every integer $d\ge2$ there are infinitely many totally real number fields of degree $d$ over which every universal quadratic form has rank at least a positive constant, depending only on $d$, times the fourth root of the discriminant. These fields are the subfields of degree $d$ of the cyclotomic fields of prime conductor $q\equiv1\pmod{2d}$.
\end{abstract}
\thanks{I would like to thank Vita Kala for his helpful comments on an earlier draft of this manuscript. The author was supported by the Czech Science Foundation GA ČR grant 26-20514S and by the Charles University programme PRIMUS/24/SCI/010.}
\maketitle

\section{Introduction}\label{sec:intro}

Let $K$ be a totally real number field with ring of integers $\OK$. A totally positive definite quadratic form over $\OK$, or more generally a quadratic $\OK$-lattice (see Section \ref{sec:prelim}), is \emph{universal} if it represents every totally positive element of $\OK$. Universal lattices exist over every totally real field by the almost local-global principle of Hsia--Kitaoka--Kneser \cite{HKK}, so
\[
m(K):=\min\{\rank\Lambda:\ (\Lambda,Q)\ \text{a universal quadratic }\OK\text{-lattice}\}
\]
is finite, and so is $m'(K)$, defined in the same way with universal quadratic forms over $\OK$, that is, free lattices. Clearly $m(K)\le m'(K)$.

That $m(K)$ can be large was first shown for real quadratic fields by Blomer--Kala \cite{BK,Ka1}, who produced for every $r$ infinitely many fields with $m(K)\ge r$. Kala--Yatsyna--\.Zmija \cite{KYZ} later proved that the real quadratic fields with $m(K)\le r$ have density zero. In degree $3$ the unboundedness is due to Yatsyna \cite{Ya} and Kala--Tinkov\'a \cite{KT}, quantitatively to Oller \cite{O}, and for multiquadratic fields to Kala--Svoboda \cite{KS} and, in a strong quantitative form, to Man \cite{Man}. Kala \cite{Ka2} then introduced a compositum technique that propagates such results from a field $L$ to $KL$ for suitable $K$, and deduced that for every $r$ there are totally real fields $F$ of degree $d$ with $m(F)\ge r$ whenever $2\mid d$ or $3\mid d$. Dole\v z\'alek \cite{Do} recast the technique group-theoretically and showed that if such fields exist in degree $d$, then they exist in every degree $kd$ with $k\ge3$. This reduces the problem to prime degrees, but for primes $d\ge5$ it remained open \cite[end of \S4]{Do}.

The method of \cite{Ya} gives lower bounds for $m(K)$ in every degree, but it requires $K$ to be primitive and to satisfy ``Condition~(A)'', namely that $K$ is monogenic and has units of every signature. For prime degrees $p\ge5$, primitivity is automatic, and Kala \cite{Ka2} noted that his compositum technique would cover all degrees if one had infinitely many totally real cyclic fields of degree $p$ containing suitable elements, for instance, via Condition~(A). However, as Kala \cite{Ka2} pointed out, M.-N. Gras \cite{Gr} proved that a cyclic field of prime degree $p\ge5$ has a power integral basis only if $2p+1$ is prime and the field is $\Q(\zeta_{2p+1})^{+}$; thus there is at most one such field in each prime degree $p\ge5$, and one must work without power integral bases or with a different Galois structure \cite{Ka2}.

We will show that Condition~(A) can be replaced by a weaker property. As already observed in \cite[\S1]{Ya}, what the method needs is a totally positive generator of the codifferent, and we will construct one explicitly for the following fields. For $d\ge2$ and a prime $q\equiv1\pmod{2d}$, let $K_q$ be the unique subfield of $\Q(\zeta_q)$ such that $[K_q:\Q]=d$. It is cyclic and totally real, with $\disc_{K_q}=q^{d-1}$ (Lemma \ref{lem:K}). For a totally real field $K$ with codifferent $\Ov$ we put
\[
M(K):=\bigl\{\beta\in\Ov:\ \beta\tp,\ \Tr\beta=1\bigr\}.
\]

\begin{theorem}\label{thm:main}
Let $d\ge2$ and let $q\equiv1\pmod{2d}$ be a prime with $q\ge d^{2}(d-1)$. Put
\[
\mu_d(q):=\frac{q^{(d-1)/2}}{2^{d-1}(d-1)!}=\frac{\disc_{K_q}^{1/2}}{2^{d-1}(d-1)!}.
\]
Then $|M(K_q)|\ge\mu_d(q)$, and:
\begin{enumerate}
\item\label{thm:main:classical} every universal classical quadratic $\mathcal O_{K_q}$-lattice has rank at least $|M(K_q)|/d\ge\mu_d(q)/d$,
\item\label{thm:main:all} if $\mu_d(q)\ge56$, every universal quadratic $\mathcal O_{K_q}$-lattice has rank at least $\dfrac{\sqrt{|M(K_q)|-56}}{d}$, which is at least $\dfrac{\sqrt{\mu_d(q)-56}}{d}$.
\end{enumerate}
\end{theorem}

Since $\disc_{K_q}=q^{d-1}$, for fixed $d$ and large $q$ we get $m'(K_q)\ge m(K_q)\gg_d\disc_{K_q}^{1/4}$, while classical lattices have rank $\gg_d\disc_{K_q}^{1/2}$.

\begin{corollary}\label{cor:everydegree}
For every $d\ge2$ and every $r$ there are infinitely many totally real number fields $K$ of degree $d$ with $m(K)\ge r$. In particular, every universal quadratic form over such $K$ has rank at least $r$.
\end{corollary}

The number of fields $K_q$ satisfying the hypotheses of Theorem~\ref{thm:main} with $\disc_{K_q} \le X$ equals the number of primes $q \equiv 1 \pmod{2d}$. By the prime number theorem in arithmetic progressions, as $X \to \infty$ this count is asymptotically
\[
\frac{d-1}{\varphi(2d)} \frac{X^{\frac{1}{d-1}}}{\ln X}.
\]
In particular, Theorem~\ref{thm:main} applies to $\gg_d X^{\frac{1}{d-1}} / \ln X$ totally real cyclic fields of degree $d$ and discriminant at most $X$.

For prime degrees $d\ge5$ this answers the question left open in \cite{Ka2,Do}. Dole\v z\'alek \cite{Do} also remarks that the compositum technique produces only Galois groups that are direct products of smaller groups, and that it would be interesting to have fields with prescribed Galois structure. The fields $K_q$ are cyclic of every degree.

\medskip

Our argument follows \cite{Ya}. If $\Ov=\delta\OK$ with $\delta\tp$ and $(\Lambda,Q)$ is a universal lattice, then $\Tr(\varepsilon\delta Q)$, with $\varepsilon=1$ for classical lattices and $\varepsilon=2$ otherwise, is a classical integral form on a $\Z$-lattice of rank $d\cdot\rank\Lambda$ with at least $2|M(K)|$ minimal vectors. An integral lattice of small rank has few minimal vectors, and this bounds $\rank\Lambda$ from below (Proposition \ref{prop:trace}). For $K_q$ we will see that $\delta=\pi/q$ for an explicit totally positive $\pi\in\OK$ (Lemma \ref{lem:K}). To show that $M(K_q)$ is large, we will identify it with the set of lattice points of a translate of $\sqrt q\,A_{d-1}$ in a simplex, and count these points by comparing the covering radius of the lattice with the inradius of the simplex (Proposition \ref{prop:count}).

\section{Preliminaries}\label{sec:prelim}

Throughout, $K$ is a totally real number field of degree $d$ with real embeddings $\sigma_1,\dots,\sigma_d$, trace $\Tr$ and norm $\Nm$. We write $\alpha\tp$ if $\alpha\in K$ is \emph{totally positive}, that is, $\sigma_i(\alpha)>0$ for all $i$, and $\OK^{+}$ for the set of totally positive elements of $\OK$. We identify $K\otimes\R$ with $\R^{d}$ by $\alpha\mapsto(\sigma_i(\alpha))_i$. Then $\OK$ is a lattice in $\R^d$ (the \emph{Minkowski embedding}), and $\Tr(\alpha\beta)=\sum_i\sigma_i(\alpha)\sigma_i(\beta)$ is the standard inner product. The \emph{codifferent}
\[
\Ov=\{\beta\in K:\ \Tr(\beta\OK)\subseteq\Z\}=\dif_K^{-1}
\]
is the dual of $\OK$ for the trace.

A \emph{quadratic $\OK$-lattice of rank $r\ge 1$} is a pair $(\Lambda,Q)$, where $\Lambda\subset K^{r}$ is a finitely generated $\OK$-module spanning $K^r$ and $Q$ is a quadratic form on $K^r$ with $Q(\Lambda)\subseteq\OK$. We write $B(u,v)=\tfrac12\bigl(Q(u+v)-Q(u)-Q(v)\bigr)$ for the associated bilinear form, so that $B(v,v)=Q(v)$. The lattice is \emph{totally positive definite} if $Q(v)\tp$ for all $v\ne0$, and \emph{universal} if it is totally positive definite and $Q(\Lambda)\supseteq\OK^{+}$. It is \emph{classical} if $B(\Lambda,\Lambda)\subseteq\OK$, that is, the bilinear form is \emph{integral}. In general, $2B(\Lambda,\Lambda)\subseteq\OK$. A quadratic form $Q=\sum_{i\le j}a_{ij}x_ix_j$ with $a_{ij}\in\OK$ is the lattice $\Lambda=\OK^r$, and it is classical if and only if $a_{ij}\in2\OK$ for $i<j$. As a $\Z$-module, $\Lambda$ is free of rank $dr$.

\section{Lower bounds from the codifferent}\label{sec:bounds}

In this section we will show that a totally positive generator of the codifferent turns a universal lattice into an integral $\Z$-lattice with many minimal vectors. We first bound the number of minimal vectors.

\begin{lemma}\label{lem:minvec}
Let $\Gamma$ be a positive definite $\Z$-lattice of rank $n\ge 1$ with an integral bilinear form $\langle\cdot,\cdot\rangle$, and let $\Phi$ be its set of non-zero vectors of minimal squared norm $\langle v,v\rangle$.
\begin{enumerate}
\item\label{lem:minvec:one} If the minimal squared norm is $1$, then $|\Phi|\le2n$.
\item\label{lem:minvec:two} If the minimal squared norm is $2$, then $|\Phi|\le2n^2+112$.
\end{enumerate}
\end{lemma}

\begin{proof}
Let $v,w\in\Phi$ with $w\ne\pm v$. By the Cauchy--Schwarz inequality, $|\langle v,w\rangle|$ is strictly smaller than the minimal squared norm.

\pr{lem:minvec:one} Here $|\langle v,w\rangle|<1$, so $\langle v,w\rangle=0$ since it is an integer. Thus $\Phi$ consists of $\pm$ pairs of pairwise orthogonal vectors, and there are at most $n$ such pairs.

\pr{lem:minvec:two} By \cite[Thm. 4.1]{RS}, a positive definite integral $\Z$-lattice of rank $n$ has at most $f(n)$ non-zero vectors of squared norm at most $2$, where $f(7)=126$, $f(n)=240+2(n-8)^2$ for $8\le n\le11$, and $f(n)=2n^2$ otherwise. In particular $|\Phi|\le f(n)\le2n^2+112$ for all $n\ge1$.
\end{proof}

\begin{proposition}\label{prop:trace}
Let $K$ be a totally real field of degree $d$ with $\Ov=\delta\OK$ for some $\delta\tp$. Then every universal $\OK$-lattice $(\Lambda,Q)$ of rank $r$ satisfies, if $|M(K)|\ge56$,
\[
r\ \ge\ \frac{\sqrt{|M(K)|-56}}{d},
\]
and if $(\Lambda,Q)$ is classical, then $r\ge|M(K)|/d$.
\end{proposition}

\begin{proof}
We may assume $M(K)$ is non-empty (otherwise $|M(K)|/d=0$, which makes the classical bound trivial, and the hypothesis $|M(K)|\ge 56$ does not hold). Put $\varepsilon=1$ if $(\Lambda,Q)$ is classical and $\varepsilon=2$ otherwise, so that $\varepsilon B(\Lambda,\Lambda)\subseteq\OK$, and let
\[
T(u,v):=\Tr\bigl(\varepsilon\delta B(u,v)\bigr)\qquad(u,v\in\Lambda).
\]
Since $\delta\in\Ov$, the form $T$ is an integral symmetric bilinear form on the $\Z$-lattice $\Lambda$ of rank $n=dr$.

For $v\in\Lambda\smallsetminus\{0\}$ the element $\delta Q(v)$ is totally positive and lies in $\delta\OK=\Ov$, so $\Tr(\delta Q(v))$ is a positive integer and $T(v,v)=\varepsilon\Tr(\delta Q(v))\ge\varepsilon$. Thus $T$ is positive definite.

Now let $\beta\in M(K)$. Then $\delta^{-1}\beta\in\OK^{+}$, so by universality there is $v_\beta\in\Lambda$ with $Q(v_\beta)=\delta^{-1}\beta$, and $T(v_\beta,v_\beta)=\varepsilon\Tr\beta=\varepsilon$. The vectors $\pm v_\beta$ for $\beta\in M(K)$ are pairwise distinct, because $Q(v_\beta)$ determines $\beta$ and $v_\beta\ne0$ (so $v_\beta\ne -v_\beta$). Hence $(\Lambda,T)$ has minimal squared norm $\varepsilon$ and at least $2|M(K)|$ minimal vectors. Lemma \ref{lem:minvec} gives $2|M(K)|\le2dr$ if $\varepsilon=1$, and $2|M(K)|\le2n^2+112$ if $\varepsilon=2$. If $\varepsilon=1$ this reads $|M(K)|\le dr$. If $\varepsilon=2$ it reads $dr\ge\sqrt{|M(K)|-56}$, since $|M(K)|\ge56$.
\end{proof}

The proposition is \cite[Prop. 5.1, Cor. 5.3, Thm. 5.4]{Ya} with the assumption that $K$ is monogenic and has units of every signature replaced by the existence of $\delta$, which is all that the argument uses \cite[\S1]{Ya}.

\section{The fields \texorpdfstring{$K_q$}{Kq}}\label{sec:fields}

Throughout this section $d\ge2$ and $q\equiv1\pmod{2d}$ is a prime. Let $\zeta=e^{2\pi i/q}$, let $H\le(\Z/q)^{\times}$ be the subgroup of index $d$, put $e=|H|=(q-1)/d$, and let $K=K_q=\Q(\zeta)^{H}$ be the fixed field of $H$. We identify $\operatorname{Gal}(\Q(\zeta)/\Q)$ with $(\Z/q)^\times$, where $a$ acts by $\zeta\mapsto\zeta^a$. In this section we will show that $\Ov$ has a totally positive generator, and that $M(K_q)$ is large.

\begin{lemma}\label{lem:K}
Let $\PP=(1-\zeta)$, let $\qq=\PP\cap\OK$, and put $\pi=\prod_{h\in H}(1-\zeta^h)$.
\begin{enumerate}
\item\label{lem:K:basic} The field $K$ is cyclic and totally real of degree $d$, with $(q)=\qq^{d}$, $\dif_K=\qq^{d-1}$ and $\disc_K=q^{d-1}$.
\item\label{lem:K:pi} The element $\pi\in\OK$ is totally positive, and $\qq=\pi\OK$.
\item\label{lem:K:codiff} The codifferent is $\Ov=q^{-1}\qq=\delta\OK$, where $\delta=\pi/q$ is totally positive, and an element $x\in\OK$ lies in $\qq$ if and only if $q\mid\Tr x$.
\item\label{lem:K:trace} The element $\eta=\sum_{h\in H}\zeta^{h}$ lies in $\OK$ and $\Tr(-\eta)=1$. In particular $\Tr(\OK)=\Z$.
\end{enumerate}
\end{lemma}

\begin{proof}
\pr{lem:K:basic} The group $(\Z/q)^\times$ is cyclic of order $q-1$ and $d\mid q-1$, so $H$ is the unique subgroup of index $d$ and $K$ is cyclic of degree $d$ over $\Q$. Since $e$ is even, $H$ contains the unique element $-1$ of order $2$, so complex conjugation fixes $K$ and $K$ is totally real.

In $\Q(\zeta)$, the prime $q$ is the only ramified prime and is totally ramified with $(q)=\PP^{q-1}$ (see \cite{Wa}). Hence $\qq$ is the unique prime ideal of $\OK$ above $q$. Ramification indices are multiplicative, $e(\PP\mid\qq)\le e$ and $e(\qq\mid q)\le d$, and their product is $q-1=de$, so $(q)=\qq^{d}$, $\qq\,\mathcal O_{\Q(\zeta)}=\PP^{e}$, and $\qq$ has residue degree $1$ and norm $q$. Since $K$ is unramified outside $q$ and $q\nmid d$, the ramification at $\qq$ is tame, so the different exponent of $\qq$ is $d-1$ \cite[Ch. III, \S6]{Se}. Hence $\dif_K=\qq^{d-1}$ and $\disc_K=\Nm(\dif_K)=q^{d-1}$.

\pr{lem:K:pi} The group $H$ permutes the factors of $\pi$, so $\pi\in K$, and $\pi$ is an algebraic integer. For $h\in(\Z/q)^\times$ the quotient $(1-\zeta^h)/(1-\zeta)$ is a unit, so each factor generates $\PP$ and $\pi\,\mathcal O_{\Q(\zeta)}=\PP^{e}=\qq\,\mathcal O_{\Q(\zeta)}$. Since an ideal of $\OK$ is determined by its extension to $\mathcal O_{\Q(\zeta)}$, $\pi\OK=\qq$.

Every real embedding of $K$ extends to an automorphism $a\in\operatorname{Gal}(\Q(\zeta)/\Q)\cong(\Z/q)^\times$, and the conjugate of $\pi$ under $a$ is $\prod_{h\in aH}(1-\zeta^{h})$. The set $aH$ is stable under $h\mapsto-h$ because $-1\in H$, and $h\ne-h$ since $q$ is odd, so $aH$ is a union of pairs $\{h,-h\}$. As
\[
(1-\zeta^{h})(1-\zeta^{-h})=|1-\zeta^{h}|^{2}>0,
\]
every conjugate of $\pi$ is positive.

\pr{lem:K:codiff} By \pr{lem:K:basic}, $\Ov=\dif_K^{-1}=\qq^{-(d-1)}=q^{-1}\qq$, since $\qq^{d}=q\OK$. By \pr{lem:K:pi} this is $(\pi/q)\OK$.

For $x\in\OK$, the inertia group, which is the whole Galois group $G$ since $\qq$ is totally ramified, acts trivially on $\OK/\qq=\mathbb F_q$ (the residue degree is $1$). Hence $\Tr(x)=\sum_{g\in G}g(x)\equiv d\,x\pmod{\qq}$. As $\Tr(x)\in\Z$ and $\qq\cap\Z=q\Z$, this gives $\Tr(x)=dx$ in $\OK/\qq=\mathbb F_q$, where $d$ is invertible because $q\nmid d$. So $x\in\qq$ if and only if $q\mid\Tr x$.

\pr{lem:K:trace} The element $\eta$ is fixed by $H$ and is an algebraic integer. Its conjugates are the sums $\sum_{h\in aH}\zeta^{h}$ for the $d$ cosets $aH$, and their sum is $\sum_{a\in(\Z/q)^\times}\zeta^{a}=-1$. So $\Tr(-\eta)=1$ and $\Tr(\OK)=\Z$.
\end{proof}

Next we recall the structure of the trace-zero lattice. Let $A_{d-1}=\{c\in\Z^{d}:\sum_ic_i=0\}$ with the standard inner product.

\begin{lemma}\label{lem:A}
The lattice $\Lambda^{0}=\{x\in\OK:\Tr x=0\}$ with the form $\Tr(xy)$ is isometric to $\sqrt q\,A_{d-1}$.
\end{lemma}
\begin{proof}
Since $\mathcal O_{\Q(\zeta)}=\Z[\zeta]$ (see \cite{Wa}) and $1+\zeta+\dots+\zeta^{q-1}=0$, the elements $\zeta,\zeta^2,\dots,\zeta^{q-1}$ span $\mathcal O_{\Q(\zeta)}$ over $\Z$. Let $C_1,\dots,C_d$ be the cosets of $H$ in $(\Z/q)^\times$, viewed as subsets of $\{1,\dots,q-1\}$. An element $x=\sum_{a=1}^{q-1}x_a\zeta^a\in\mathcal O_{\Q(\zeta)}$ ($x_a\in\Z$) lies in $\OK=\mathcal O_{\Q(\zeta)}\cap K=\mathcal O_{\Q(\zeta)}^H$ if and only if
\[
\sum_{a=1}^{q-1}x_a\zeta^{ha}=\sum_{a=1}^{q-1}x_a\zeta^a\qquad\text{for all }h\in H,
\]
which means that $a\mapsto x_a$ is constant on each coset $C_i$. Hence
\[
\eta_i:=\sum_{a\in C_i}\zeta^a\qquad(i=1,\dots,d)
\]
form a $\Z$-basis of $\OK$.

Since $\operatorname{Gal}(K/\Q)\cong(\Z/q)^\times/H$ permutes the cosets $C_1,\dots,C_d$ simply transitively, each $\eta_i$ has trace
\[
\Tr(\eta_i)=\sum_{j=1}^d\eta_j=\sum_{a=1}^{q-1}\zeta^a=-1.
\]
Thus an element $x=\sum_{i=1}^d c_i\eta_i\in\OK$ ($c_i\in\Z$) satisfies $\Tr x=-\sum_{i=1}^d c_i$, so $x\in\Lambda^0$ if and only if $c=(c_1,\dots,c_d)\in A_{d-1}$.

To compute $\Tr(\eta_i\eta_j)$, recall that for any $m\in\Z$,
\[
\Tr_{\Q(\zeta)/\Q}(\zeta^m)=\sum_{a=1}^{q-1}\zeta^{am}=
\begin{cases}
q-1 & \text{if } q\mid m,\\
-1 & \text{if } q\nmid m.
\end{cases}
\]
Because $e=[\Q(\zeta):K]$ and $-1\in H$, we have $-C_i=C_i$ for every $i$. In particular, for $(a,b)\in C_i\times C_j$ we have $q\mid(a+b)$ if and only if $b\equiv -a\pmod q$, which occurs for $e$ pairs $(a,-a)$ when $i=j$, and for no pairs when $i\ne j$. Since $\Tr_{\Q(\zeta)/\Q}(\alpha)=e\Tr(\alpha)$ for $\alpha\in K$, it follows that
\[
e\Tr(\eta_i^2)=\Tr_{\Q(\zeta)/\Q}(\eta_i^2)=e(q-1)+(e^2-e)(-1)=e(q-e),
\]
and for $i\ne j$,
\[
e\Tr(\eta_i\eta_j)=\Tr_{\Q(\zeta)/\Q}(\eta_i\eta_j)=e^2(-1)=-e^2.
\]
Dividing by $e$, we obtain $\Tr(\eta_i^2)=q-e$ and $\Tr(\eta_i\eta_j)=-e$ for $i\ne j$.

Therefore, for any $x=\sum_{i=1}^d c_i\eta_i$ and $y=\sum_{j=1}^d c'_j\eta_j$ in $\Lambda^0$ (so $c,c'\in A_{d-1}$, meaning $\sum_i c_i=\sum_j c'_j=0$),
\[
\Tr(xy)
=\sum_{i,j=1}^d c_ic'_j\Tr(\eta_i\eta_j)
=q\sum_{i=1}^d c_ic'_i-e\Bigl(\sum_{i=1}^d c_i\Bigr)\Bigl(\sum_{j=1}^d c'_j\Bigr)
=q\sum_{i=1}^d c_ic'_i.
\]
Consequently, $\sqrt{q}\,c\mapsto\sum_{i=1}^d c_i\eta_i$ is an isometry from $\sqrt{q}\,A_{d-1}$ onto $(\Lambda^0,\Tr(xy))$.
\end{proof}
For prime $d$ this also appears in \cite[Thm. 2.9]{MM} (see also Conner--Perlis \cite{CP}), and for general $d$ see also \cite[Thm. 2.10]{BMS2}.

For the count, we need the following lemma on lattice points in convex sets.

\begin{lemma}\label{lem:convex}
Let $E$ be a Euclidean affine space of dimension $n\ge1$, let $\Lambda$ be a lattice of rank $n$ in its space of translations with covering radius $\rho$, and let $L\subset E$ be a coset of $\Lambda$. Let $C\subset E$ be a compact convex set containing a closed ball of radius $R>\rho$, and put $\tau=\rho/R$. Then
\[
(1-\tau)^{n}\,\frac{\vol(C)}{\covol(\Lambda)}\ \le\ |L\cap C|\ \le\ (1+\tau)^{n}\,\frac{\vol(C)}{\covol(\Lambda)}.
\]
\end{lemma}

\begin{proof}
We will show that the union of the Voronoi cells around the points of $L\cap C$ lies between a slightly smaller and a slightly larger copy of $C$, both homothetic to $C$ about the centre of a ball in $C$.

Let $c$ be the centre of a ball $\mathcal B(c,R)\subseteq C$. Let $V$ be the Voronoi cell of $\Lambda$, that is, the set of points that are at least as close to $0$ as to every other point of $\Lambda$. We recall that $V=-V$ is a convex polytope (so its boundary has measure zero), that the translates $x+V$ with $x\in L$ cover $E$ and have pairwise disjoint interiors, that $\vol(V)=\covol(\Lambda)$, and that $\|z\|\le\rho$ for $z\in V$. Also $L\cap C$ is finite because $L$ is discrete and $C$ is compact, and $\rho>0$ (since $n\ge1$), so $0<\tau<1$. For $s>0$ let $C_s=c+s(C-c)$, so that $\vol(C_s)=s^{n}\vol(C)$.

Let $z$ be a vector with $\|z\|\le\rho$. Then $z=\tau w$ with $\|w\|\le R$, so $c+w\in C$. We will use this in both inequalities.

For the lower bound, let $y\in C_{1-\tau}$, so $y=(1-\tau)p+\tau c$ for some $p\in C$. The point $y$ lies in a cell $x+V$ with $x\in L$, so $x-y\in -V=V$ and $x=y+\tau w$ with $\|w\|\le R$. Then
\[
x=(1-\tau)p+\tau(c+w)\in C
\]
by convexity. So the cells $x+V$ with $x\in L\cap C$ cover $C_{1-\tau}$, and $(1-\tau)^{n}\vol(C)\le|L\cap C|\covol(\Lambda)$.

For the upper bound, let $x\in L\cap C$ and $y\in x+V$, so $y=x+\tau w$ with $\|w\|\le R$. Then
\[
c+\frac{y-c}{1+\tau}=\frac{1}{1+\tau}\,x+\frac{\tau}{1+\tau}\,(c+w)\in C,
\]
so $y\in C_{1+\tau}$. The cells $x+V$ with $x\in L\cap C$ have pairwise disjoint interiors and lie in $C_{1+\tau}$, so $|L\cap C|\covol(\Lambda)\le(1+\tau)^{n}\vol(C)$.
\end{proof}

\begin{proposition}\label{prop:count}
Assume $q\ge d^{2}(d-1)$ and put $\theta=d\sqrt{d-1}/(2\sqrt q)$. Then $\theta\le\frac12$ and
\[
(1-\theta)^{d-1}\,\frac{q^{(d-1)/2}}{(d-1)!}\ \le\ |M(K_q)|\ \le\ (1+\theta)^{d-1}\,\frac{q^{(d-1)/2}}{(d-1)!}.
\]
\end{proposition}

\begin{proof}
We will identify $M(K_q)$ with the set of points of a translate of $\Lambda^0$ in a simplex, and then apply Lemma \ref{lem:convex}.

We have $\theta\le\frac12$ if and only if $q\ge d^{2}(d-1)$. Let $\mathbf e_1,\dots,\mathbf e_d$ be the standard basis and $\mathbf 1=(1,\dots,1)$. Let $\mathcal H_0=\{y\in\R^d:\sum_jy_j=0\}$ and $\mathcal H=\{y\in\R^d:\sum_jy_j=q\}$, and let
\[
\mathcal S=\{y\in\mathcal H:\ y_j>0\text{ for all }j\},\qquad L=\{x\in\OK:\ \Tr x=q\}\subset\mathcal H
\]
under the Minkowski embedding. We first show that $|M(K)|=|L\cap\mathcal S|$ and that $L$ is a coset of $\Lambda^0$.

If $\beta\in M(K)$, then $x=q\beta$ lies in $\qq\subseteq\OK$ by Lemma \ref{lem:K}\pr{lem:K:codiff}, and it is totally positive with $\Tr x=q$, so $x\in L\cap\mathcal S$. Conversely, if $x\in L\cap\mathcal S$, then $q\mid\Tr x$, so $x\in\qq$ by Lemma \ref{lem:K}\pr{lem:K:codiff}, and $x/q\in\Ov$ is totally positive of trace $1$. Thus $\beta\mapsto q\beta$ is a bijection $M(K)\to L\cap\mathcal S$. Moreover $x_0=-q\eta\in\OK$ has $\Tr x_0=q$ by Lemma \ref{lem:K}\pr{lem:K:trace}, so $L$ is non-empty, and hence $L=x_0+\Lambda^{0}$. Finally $0\notin L$, and a non-zero element of $\OK$ has no zero conjugate, so $L$ does not meet the boundary of $\overline{\mathcal S}=\{y\in\mathcal H:y_j\ge0\text{ for all }j\}$, and $L\cap\overline{\mathcal S}=L\cap\mathcal S$.

\medskip

We now apply Lemma \ref{lem:convex} to $E=\mathcal H$, $n=d-1$, the lattice $\Lambda^{0}\subset\mathcal H_0$, the coset $L$ and $C=\overline{\mathcal S}$. Let $\widehat y=(q/d,\dots,q/d)\in\mathcal H$ and $R=q/\sqrt{d(d-1)}$. If $y\in\mathcal H$ and $\|y-\widehat y\|\le R$, then $y-\widehat y$ is orthogonal to $\mathbf 1=(1,\dots,1)$, so
\[
y_j=\frac qd+\bigl\langle\mathbf e_j-\tfrac1d\mathbf 1,\ y-\widehat y \bigr\rangle\ \ge\ \frac qd-\sqrt{\frac{d-1}{d}}\,R=0,
\]
using $\|\mathbf e_j-\frac1d\mathbf 1\|^{2}=(d-1)/d$. This shows that $\overline{\mathcal S}$ contains the ball of radius $R$ about $\widehat y$ in $\mathcal H$.

Under the Minkowski embedding, $\Tr(xy)$ is the standard Euclidean inner product on $\R^d$. By Lemma \ref{lem:A}, $\Lambda^0\cong\sqrt q\,A_{d-1}$, so its covering radius is $\rho=\sqrt q\,\rho(A_{d-1})$. By \cite[Ch. 4, \S6.1]{CS}, $\rho(A_{n})^{2}=a(n+1-a)/(n+1)$ with $a=\lfloor(n+1)/2\rfloor$, and $a(n+1-a)\le(n+1)^{2}/4$. Hence $\rho(A_{d-1})\le\sqrt d/2$ and
\[
\tau=\frac\rho R\ \le\ \frac{\sqrt{qd}}{2}\cdot\frac{\sqrt{d(d-1)}}{q}=\theta\ \le\ \frac12 .
\]
Thus $R>\rho$, and Lemma \ref{lem:convex}, together with $\tau\le\theta$ and the monotonicity of $(1\pm t)^{d-1}$ in $t$, gives
\[
(1-\theta)^{d-1}\,\frac{\vol(\overline{\mathcal S})}{\covol(\Lambda^0)}\ \le\ |L\cap\overline{\mathcal S}|\ \le\ (1+\theta)^{d-1}\,\frac{\vol(\overline{\mathcal S})}{\covol(\Lambda^0)} .
\]

It remains to compute the volumes. Under the orthogonal projection forgetting $y_d$, the set $\overline{\mathcal S}$ projects onto $\{y\in\R^{d-1}:\ y_j\ge0\text{ for all }j,\ \sum_{j<d}y_j\le q\}$. This simplex has volume $q^{d-1}/(d-1)!$, and the projection multiplies $(d-1)$-dimensional volumes in $\mathcal H$ by the cosine $1/\sqrt d$ of the angle between the normals $\mathbf 1/\sqrt d$ and $\mathbf e_d$. Hence $\vol(\overline{\mathcal S})=\sqrt d\,q^{d-1}/(d-1)!$. Also, $\covol(A_{d-1})=\sqrt d$ \cite[Ch. 4, \S6.1]{CS}, so $\covol(\Lambda^0)=q^{(d-1)/2}\covol(A_{d-1})=\sqrt d\,q^{(d-1)/2}$, and the volume ratio is $q^{(d-1)/2}/(d-1)!$.
\end{proof}

Since $\theta\le\frac12$, the proposition gives $|M(K_q)|\ge\mu_d(q)>0$, and $|M(K_q)|\sim q^{(d-1)/2}/(d-1)!$ for fixed $d$ as $q\to\infty$.

\section{Proofs of the main results}\label{sec:proofs}

\begin{proof}[Proof of Theorem \ref{thm:main}]
By Lemma \ref{lem:K}\pr{lem:K:basic}, $K_q$ is totally real of degree $d$ with $\disc_{K_q}=q^{d-1}$, and by Lemma \ref{lem:K}\pr{lem:K:codiff}, $\Ov=\delta\OK$ with $\delta\tp$. Proposition \ref{prop:count} and $\theta\le\frac12$ give
\[
|M(K_q)|\ \ge\ (1-\theta)^{d-1}\frac{q^{(d-1)/2}}{(d-1)!}\ \ge\ \frac{q^{(d-1)/2}}{2^{d-1}(d-1)!}=\mu_d(q).
\]
Now Proposition \ref{prop:trace} gives \pr{thm:main:classical} and \pr{thm:main:all}, and the second inequalities follow from $|M(K_q)|\ge\mu_d(q)$.
\end{proof}

\begin{proof}[Proof of Corollary \ref{cor:everydegree}]
Fix $d\ge2$ and $r$. By Dirichlet's theorem there are infinitely many primes $q\equiv1\pmod{2d}$, and $\mu_d(q)\to\infty$ as $q\to\infty$. For all such $q$ with $q\ge d^{2}(d-1)$, $\mu_d(q)\ge56$, and $\bigl(\sqrt{\mu_d(q)-56}\bigr)/d\ge r$, Theorem \ref{thm:main}\pr{thm:main:all} gives $m(K_q)\ge r$. These fields are pairwise distinct, since $\disc_{K_q}=q^{d-1}$. As $m(K)\le m'(K)$, every universal quadratic form over $K_q$ has rank at least $r$.
\end{proof}

\end{document}